\documentclass[12pt,a4paper]{article}
\usepackage[a4paper,margin=1in]{geometry}
\usepackage{amsmath,amssymb}
\usepackage{graphicx}
\usepackage{setspace}
\usepackage{tikz}
\usetikzlibrary{arrows.meta,positioning,matrix}
\usepackage[hidelinks]{hyperref}
\usepackage{xurl}

\newcommand{\citetext}[1]{(#1)}
\newenvironment{reflist}{%
  \begingroup\setstretch{1.1}\setlength{\parindent}{-1.5em}%
  \setlength{\leftskip}{1.5em}\setlength{\parskip}{6pt}\small}{\endgroup}

\begin{document}

\begin{center}
  {\LARGE\bfseries Ascertaining Learning:}\\[6pt]
  {\large\bfseries Generative AI, the Proxy Problem, and Cognitive Competence in University Mathematics}\\[10pt]
  {\large B.~Osano}\\[4pt]
  {\small Centre for Higher Education Development \&\\
  Cosmology and Gravity Group, Department of Mathematics and Applied Mathematics\\
  University of Cape Town, Private Bag X3, Rondebosch 7701, South Africa\\
  \href{mailto:bob.osano@uct.ac.za}{bob.osano@uct.ac.za}}
\end{center}
\vspace{12pt}

\begin{abstract}
\noindent
For decades, university mathematics assessed \textit{executable procedure}, computing and
reproducing standard proofs, as a proxy for the understanding it cared about, because procedure
was cheap to mark and understanding was not. Generative AI now executes much of this
procedure. It does not merely let students bypass the proxy; it \textit{exposes}
how much of what was certified was a proxy. The crisis is not simply cheating but a \textit{proxy
problem}, and resolving it forces a question the discipline has long avoided: what is the
irreducible human core of mathematical competence, and how can it be built when some of the
executable work can be offloaded? We set out two futures: \textit{displacement-and-re-centering},
in which AI takes over some execution and education moves toward judgment, formulation and
verification; and \textit{atrophy}, in which premature offloading may erode capacities that depend
on prior construction. At the heart of the optimistic future is a \textit{verification instability}:
some forms of oversight are cheap and bounded, while others approach the competence required to
construct the work. The outcome is a design choice, not a technological destiny. Rather than add a model,
we state the response as two principles, \textit{construct before you delegate}, and \textit{let
verification set the ceiling}: how far a student can independently verify AI output limits what they
may safely delegate. We illustrate with a sequence and two worked problems, and close on the
question: whether the construct-before-delegate dependency is a hard constraint or an artefact of
how mathematics has been taught.

\end{abstract}

\section{Introduction: A Discipline Discovers Its Proxy}
Assessment in undergraduate mathematics is guided, largely unintentionally, by competing philosophies of what mathematics is. Social constructivism \citetext{Ernest, 1991, 1998} treats mathematical knowledge as a socially negotiated activity, shifting assessment toward communicating arguments and constructing meaning. Formalism \citetext{Curry, 1951; Shapiro, 2000} treats it as a system of rules and rigorous deductions, guiding exams that test proof-writing and logical validity. Pragmatism \citetext{Ernest, 1989; Thompson, 1992} treats it as a tool, shaping applied assessments of modelling for real-world problems. Together these dictate whether tests evaluate algorithmic, procedural competence or conceptual, proof-oriented understanding.

Consider a standard first-year calculus problem: evaluate a definite integral by substitution and integration by parts. For most of the twentieth century, a student who produced the correct answer with correct working had shown the examiner something real, not the answer, which is of marginal interest, but the understanding and procedural fluency it stood \textit{for}. The mark was never about the integral. It was about the competence the integral was a \textit{proxy} for, what measurement theory calls the gap between an observable performance and the latent construct it is meant to indicate \citetext{Messick, 1989}. And the proxy was used because it was cheap. Procedure can be set, sat, and marked at scale; understanding, judgment, and reasoning, which carry long-term learning \citetext{Boud \& Falchikov, 2006}, are expensive to elicit and expensive to grade. Mathematics education, like most of education, built its assessment machinery on proxies of convenience, in frequent breach of the constructive-alignment principle that tasks should elicit the outcomes they claim to measure \citetext{Biggs, 1996}. So the question in the title is not idle. When an examination like this is set, what are we actually \textit{ascertaining}: that a student has learned, or only that a correct answer has appeared?

Generative AI can destabilise the proxy. A system that performs the substitution, the integration
by parts, and the layout in seconds can sever the evidential link between the artefact (correct
working) and the competence it was intended to indicate. The first institutional reaction has often
been to treat this as cheating, students using a forbidden tool, to be met with prohibition and
detection. That response addresses tool use, but it can miss the deeper measurement problem. The
issue is that the proxy has been revealed as a proxy: once a machine can contribute substantially to
the observed performance, the performance no longer warrants the same inference about the learner's
competence without additional evidence.

The reframing is clarifying, not merely destructive. The question AI poses is not ``how do we stop
students using it?'' but an older and deeper one the discipline could postpone only while the
proxy worked: \textit{what is the irreducible human core of mathematical competence, and how do we
build and certify it?} As long as executable procedure was an adequate stand-in, mathematics
education did not have to answer sharply. It can no longer avoid it.

The measurement problem can be stated more explicitly as a validity chain:
\[
\text{construct} \;\longrightarrow\; \text{assessment task} \;\longrightarrow\;
\text{observable performance} \;\longrightarrow\; \text{inference about competence}.
\]
Under unaided assessment, the final inference holds only insofar as performance indicates the
construct. With AI, the observed performance may be produced jointly by learner competence and
machine assistance, so the pathway is better written as
\[
\text{student competence}+\text{AI assistance}
\;\longrightarrow\;\text{observed performance}.
\]
This does not make performance meaningless, but it changes what can validly be inferred: the proxy
problem is one of evidential warrant for the inference, not a claim that every AI-assisted artefact is
invalid.

\subsection{Related Work and the Gap}

Four separate literatures bear on the problem, and each covers only part of it. Work on
\textit{AI capability in mathematics} documents what generative and neural systems can now do
\citetext{Frieder et al., 2023; Trinh et al., 2024; AlphaProof, 2024} but is silent on the
educational response. A literature on \textit{AI-era reframings of cognitive taxonomies} updates
Bloom's taxonomy for AI \citetext{Hmoud \& Shaqour, 2024; Ayodele et al., 2026; Jain \& Samuel,
2025; Gonsalves, 2026; Clark, 2025}, but redefines what students do \textit{with} AI at each level
while assuming the lower, executable levels can be safely handed to the machine, sidestepping
the developmental question of whether the higher levels survive that delegation, and never treating
verification as a distinct competence. The \textit{cognitive science of learning}, cognitive
offloading \citetext{Risko \& Gilbert, 2016; Sparrow et al., 2011}, productive failure and
desirable difficulties \citetext{Bjork \& Bjork, 2011; Kapur, 2008}, and cognitive load
\citetext{Sweller, 1988}, establishes how removing effortful construction damages durable
learning, but predates the generative-AI shift and supplies the mechanism of harm without a design
response. And work on \textit{assessment and the human factors of automation}, validity theory
\citetext{Messick, 1989}, constructive alignment \citetext{Biggs, 1996}, and the ``ironies of
automation'' \citetext{Bainbridge, 1983; Parasuraman \& Manzey, 2010}, shows how assessment can
measure a proxy and how automating routine work erodes the skill needed to oversee it, but never
brings that oversight problem to bear on \textit{building} mathematical competence.

The gap lies at their common intersection. No existing account (i)~diagnoses the disruption as the
exposure of assessment's proxy problem rather than as cheating; (ii)~shows that the popular
``re-centre on higher-order skills'' response collides with a \textit{developmental-dependency
bottleneck} and a \textit{verification instability}; or (iii)~offers a principled way to sequence
human construction before AI delegation and to gate delegation on verification competence. This
paper makes those three moves, states the design response as two principles rather than another
named framework, and names, rather than conceals, the field's central open question: whether
the construct-before-delegate dependency is a hard constraint or a contingent artefact of past
teaching.

\section{How the Taxonomies Encoded the Proxy}

\subsection{Bloom and the Machinery of Assessment}
Bloom's revised taxonomy \citetext{Anderson \& Krathwohl, 2001} has been widely used to organise
curriculum and assessment design, including in mathematics. Bloom orders cognition from
\textit{Remember} and \textit{Understand} through \textit{Apply} and \textit{Analyse} to
\textit{Evaluate} and \textit{Create}. It was built to help educators specify what students should be
able to do, and to design assessment that elicited it.

There is, however, a difficulty when these general taxonomies are treated as accounts of mathematical
thinking. Long, Dunne, and de Kock \citetext{2014} found that Bloom's categories can provide useful
general descriptions but are difficult to interpret meaningfully in mathematics, because the
cognitive demand of a task depends on its mathematical content and on learners' prior experience.
Darlington \citetext{2013}, writing specifically about advanced mathematics questions, likewise
notes the longstanding criticism that Bloom is an awkward fit for mathematics and contrasts it with
the mathematics-specific MATH taxonomy. This is not an argument that Bloom has no value. Indeed,
Radmehr and Drake \citetext{2019} argue that Revised Bloom's Taxonomy can be useful in mathematics,
particularly because its two-dimensional structure incorporates metacognitive knowledge and does
not require a strict hierarchy. The more limited point is that a Bloom category does not, by itself,
tell us what mathematical competence a student possesses.

The distinction matters because assessment economies can turn a descriptive taxonomy into a proxy
for competence. The lower and middle levels, remembering definitions, understanding concepts,
applying procedures, are cheaper to assess at scale than the judgment and reasoning they are meant to
support, so assessment gravitates toward their most tractable evidence: correct execution of standard
procedure. The empirical record bears this out: Bergqvist \citetext{2007} found most examined tasks
solvable by \textit{imitative} rather than creative reasoning, and Lithner \citetext{2008} showed how
routinely assessment rewards reproduced procedure over reasoning proper. A curriculum nominally aimed
at \textit{Create} can, in what is marked, become one that certifies \textit{Apply}. The taxonomies
did not cause this; the economics of assessment did. But they supplied a vocabulary in which the
observable activity could pass as evidence of the intended construct.

Mathematics education has, in fact, produced a taxonomy that makes this distinction more explicit.
Smith and colleagues' \citetext{Smith et al., 1996} Mathematical Assessment Task Hierarchy, the
MATH taxonomy, was developed specifically to help construct examinations that assess a range of
mathematical knowledge and skills. It distinguishes Group~A (factual recall, comprehension, and
routine use of procedures), Group~B (information transfer and application in new situations), and
Group~C (justification, interpretation, conjecture, and evaluation). Its importance here is not
that it provides a replacement hierarchy for Bloom, but that it classifies mathematical tasks by
the \textit{nature of the activity} rather than assuming that a generic cognitive level captures
mathematical demand. Darlington \citetext{2013} similarly uses MATH to contrast the kinds of
questions posed in advanced mathematics examinations. Long et al. \citetext{2014} go further in
arguing that mathematical dimensions need not be hierarchical: within a particular mathematical
dimension, a task may require anything from recall to highly creative activity. The point is
therefore not that mathematics needs a better ladder of cognitive levels. It is that mathematical
activity has structure that a generic ladder cannot, by itself, represent.

Read this way, the MATH groups express a dependency of \textit{competence} rather than a fixed order
of performance: transfer and application (Group~B) presuppose the automated knowledge and routines of
Group~A, and justification, interpretation, and evaluation (Group~C) presuppose Group~B. This
captures the developmental dependency without imposing the linear task order that mathematics does
not honour. A student may begin a problem by formulating it, an act high on any generic ladder, and
still be unable to carry it through without the automated lower-order competence beneath it. The order
in which mathematics is done can be non-linear and recursive while the dependency between kinds of
competence still holds.

The proxy operates even where assessment believes it has climbed above procedure. Students routinely
learn to \textit{state} the $\varepsilon$--$\delta$ definition of a limit yet cannot \textit{use} it,
structure a proof from it, or decide whether a function satisfies it \citetext{Tall \& Vinner, 1981}.
The gap between artefact and construct then opens \textit{inside} what looks like conceptual
assessment, so moving assessment ``up'' the range does not, by itself, escape the proxy.

The distinction is important for the argument that follows. A mathematics-specific task
classification can improve the design of assessment, and a revised Bloom taxonomy can provide a
useful vocabulary for analysing the cognitive activity an assessment invites. Neither, however,
establishes that the performance observed is evidence of the underlying mathematical competence.
The paper's concern is therefore not to replace Bloom with another taxonomy. It is to separate the
description of a task from the inference made from performance on that task. Generative AI makes
that separation newly urgent because it can now produce much of the observable performance without
the learner necessarily possessing the competence the performance was intended to evidence.

\subsection{Why This Matters Now}

The relevance to generative AI is direct. AI is highly capable exactly where assessment has focused,
the executable procedures that have served as convenient proxies, so the disruption falls hardest on
performance that is cheap to specify and assess. The taxonomies remain useful for re-weighting the
activities they describe, but they cannot by themselves say what to do, because classifying the
activity a task invites is not evidence that the learner possesses the competence the task was meant
to reveal. That distinction is the question the rest of the paper takes up.

\section{What Generative AI Can Now Do}

\subsection{A Capability Audit}

Today's large language and reasoning models, together with symbolic engines and proof assistants,
can attempt tasks across a wide range of undergraduate mathematical activity, but their reliability
is highly task-dependent. They are often strong on recall and standard procedure, while systematic
benchmarking still finds meaningful error rates on non-trivial mathematics \citetext{OpenAI, 2023;
Frieder et al., 2023}. They have also advanced substantially on some forms of mathematical proof:
AlphaGeometry solved 25 of 30 selected olympiad-level geometry problems, approaching the performance
of an average IMO gold medallist \citetext{Trinh et al., 2024}. These results do not establish a
simple ordering of AI capability by cognitive level. Rather, they show uneven \textit{delegability}:
some mathematical activities can be delegated with relatively clear checks, while others involve
assumptions, interpretation, or standards of correctness that are much harder to externalise. The
educational question is therefore not where AI sits on a cognitive ladder, but which components of a
mathematical task can be safely delegated and what competence the learner must retain to evaluate the
result.

The capability literature therefore supports a narrower claim: AI can now contribute to mathematical
activity well beyond routine calculation, but unevenly and task-specifically. That is all the argument
needs. It requires no claim about a single frontier of what AI can and cannot do, only that the
machine can supply enough of the observable work that assessment can no longer assume the learner
alone produced it.

\subsection{The Uneven Shape of the Disruption}

The audit's key finding is therefore not a single capability level or a Bloom-shaped profile, but
an uneven distribution of tasks that can be delegated (Figure~1). Some tasks are readily delegated because
their outputs are structured and their correctness can be checked locally. Others remain difficult to
delegate safely because correctness depends on modelling assumptions, interpretation, or open-ended
judgment. The central tension of the paper arises when the tasks most easily delegated are also tasks
whose prior construction contributes to the competence needed to understand or verify later work.

\begin{center}
\vspace{4pt}
\begin{tikzpicture}[>=Latex]
  \node[draw, rounded corners=2pt, align=center, text width=2.9cm, inner sep=4pt, minimum height=1.3cm]
    (bounded) at (0,0) {\textbf{Structured task}\\[2pt]{\footnotesize clear output;\\ bounded check}};
  \node[draw, rounded corners=2pt, align=center, text width=2.9cm, inner sep=4pt, minimum height=1.3cm]
    (mixed) at (5.6,0) {\textbf{Partly open task}\\[2pt]{\footnotesize multiple methods;\\ partial checks}};
  \node[draw, rounded corners=2pt, align=center, text width=2.9cm, inner sep=4pt, minimum height=1.3cm]
    (open) at (11.2,0) {\textbf{Open-ended task}\\[2pt]{\footnotesize model choice;\\ interpretation}};
  \draw[->, line width=0.8pt] (bounded.east) -- (mixed.west)
    node[midway, above, font=\footnotesize]{less bounded};
  \draw[->, line width=0.8pt] (mixed.east) -- (open.west)
    node[midway, above, font=\footnotesize]{more judgment};
  \node[align=center, font=\footnotesize, text width=3.2cm] at (0,-1.9)
    {AI can often supply\\ mechanical sub-tasks};
  \node[align=center, font=\footnotesize, text width=3.2cm] at (5.6,-1.9)
    {delegation requires\\ stronger oversight};
  \node[align=center, font=\footnotesize, text width=3.2cm] at (11.2,-1.9)
    {delegation may require\\ near-full competence};
\end{tikzpicture}

\vspace{2pt}
{\small\textbf{Figure 1: Task-dependent delegability.} The figure does not rank AI by cognitive
level. It illustrates the relevant distinction for this paper: the safety of delegation depends on
the structure of the task and, especially, on whether the learner has a bounded way to verify the
output. As tasks become more open-ended, verification increasingly depends on modelling judgment,
interpretation, and domain knowledge.\par}
\vspace{4pt}
\end{center}

\section{\label{sec4}Two Futures: Displacement-and-Re-centering versus Atrophy}

The disruption admits two coherent readings. They point in opposite directions, and between them
they define the field's decision.

\subsection{Displacement-and-Re-centering}

On the optimistic reading, AI takes over the routine execution layer and frees mathematics
education to move up the range, to the capacities AI is weak at and practitioners still need:
problem formulation and modelling (turning an ill-defined situation into tractable mathematics,
which needs an understanding of the world as well as the mathematics); mathematical judgment and
taste (which method to try, whether an answer is even plausible); the \textit{comprehension and
critique} of proof rather than its mechanical production; verification and error detection;
transfer to novel contexts; and the communal construction and communication of understanding. On
this reading the outlook is favourable. The discipline stops drilling what machines do better and
re-centres on the irreducibly human. There is precedent for this. The calculator transition
re-centred numeracy education on number sense rather than hand computation, and the feared collapse
of competence did not follow \citetext{Ellington, 2003}.

\subsection{Atrophy}

On the pessimistic reading, the optimistic list is developmentally incoherent. The higher
capacities it names are not free-standing; they are built \textit{on} the executable base. A
student cannot exercise judgment about an integral they could never compute, cannot critique a
proof they could not have constructed, cannot spot an error in a derivation whose correct form they
never internalised. Remove the base, let students offload execution before the schemas that
execution builds have formed, and the higher capacities have nothing to stand on. On this
reading, ``let AI do the routine and teach humans the judgment'' is not a re-centring but a loss of
the underlying competence: graduates fluent with AI and incompetent without it, whose credential
certifies a competence the human no longer holds.

\subsection{Why the Calculator Precedent Does Not Settle It}

The optimistic reading leans on a real precedent, and it deserves a head-on answer. Arithmetic
calculators, and later computer-algebra systems, provoked the same fear of competence collapse and
did not, on balance, produce it: numeracy education re-specified itself around number sense
\citetext{Ellington, 2003}. But read carefully, that history does not refute the argument of this
paper; it illustrates it. Offloading stayed safe under two conditions: the base was constructed
first, and the tool was confined to a bounded layer.

Even the clearest loss fits this pattern. Mental arithmetic did suffer, but that cost is not uniform
across a learner's development. At the foundational stage, mental computation is how learners first
build number sense and learn to reason without external aid, and offloading it there would remove the
very work that constructs the base. Beyond that stage, once the base is in place, the same mechanical
skill can be offloaded without loss of the conceptual understanding where competence lies. University
students sit well past that point, so a calculator carrying the mechanical load costs them nothing
conceptual. Integration survived for the same reason: the teaching set the base before the tool was
allowed to carry it.

Generative AI keeps the first condition but breaks the second, on two features. The first is
\textit{scope}: a calculator has a hard ceiling, it cannot formulate a model, construct a proof, or
explain a result, so the layers above mechanical computation still had to be built by hand, and
re-specification worked precisely because the tool's incapacity protected them. Generative AI offers
plausible output across the \textit{entire} range, so no layer is protected by the tool's own limits.
The second is \textit{legibility}: a calculator is reliable within its domain and silent outside it,
whereas generative AI is fluently, confidently wrong in ways a novice cannot detect and does not mark
the boundary of its own competence \citetext{Ji et al., 2023}, and people under-scrutinise such
output, treating fluency or the mere presence of an explanation as a signal of competence
\citetext{Bu\c{c}inca, Malaya, \& Gajos, 2021}. The step-by-step explanation aids learning when it is
correct and is exactly what makes an error persuasive when it is not.

None of this proves the pessimistic reading. It shows only that the analogy cannot do the reassuring
work asked of it: the features that made the calculator transition safe, protected scope and legible
reliability, are the ones generative AI lacks. Whether delegation is nonetheless safe has to be
earned elsewhere, from whether the base is constructed first and whether the student can independently
verify the output, the empirical question Section~5 isolates. The precedent cannot be invoked as
though the answer were known.

\section{\label{sec5}The Developmental-Dependency Bottleneck}

\subsection{The Constructive Character of Mathematical Knowledge}

The pessimistic reading rests on a well-supported account of how mathematical competence is
acquired. A broad tradition of learning theory holds that abstract mathematical objects are
\textit{constructed} through activity on their more concrete precursors, Vygotsky's
\citetext{1978} account of concept formation through socially mediated activity and Sfard's
\citetext{1991} process--object duality being two influential statements. In Sfard's terms, a mathematical object begins life as a process the learner performs and becomes an
object of thought only once that performance is interiorised: the concept of a function is built by
working with functions, the idea of a limit by computing limits, a proof technique by executing it,
failing, and executing it again. Dubinsky's \citetext{1991} reflective abstraction and Tall and
Vinner's \citetext{1981} concept image describe the same dependence of understanding on prior
activity. Skip the performance and the object does not form. Cognitive load theory
\citetext{Sweller, 1988; Kirschner et al., 2006} tells the same story from the architecture of
memory: schemas are built only through the effortful processing that routine offloading removes.
And that this construction cannot be shortcut without cost is the convergent finding of research on
productive struggle and productive failure \citetext{Hiebert \& Grouws, 2007; Kapur, 2008} and on
desirable difficulties \citetext{Bjork \& Bjork, 2011}, whose unsettling result is that the
conditions making learning feel easiest are often those that build the least durable competence.

\subsection{The Bottleneck}

The consequence is a possible bottleneck. Higher-order mathematical competence may depend
substantially on the conceptual and procedural schemas developed through prior engagement with the
underlying mathematics. This is not the claim that execution is sufficient for understanding; it
plainly is not, nor is it the claim that every mathematical procedure must first be mastered by
hand. The more defensible claim is that students need sufficiently developed domain-specific
representations and procedures to understand, judge, and verify later mathematical work, and that
historically much of this development has occurred through active procedural engagement. Convergent
causal evidence points the same way: prior mathematical knowledge exerts a causal effect on later
performance, mediated by self-efficacy, in introductory undergraduate courses
\citetext{Zakariya, 2021}. Whether AI
can support the same developmental functions through other routes is an open empirical question.
Everything therefore turns on whether the dependency is a \textit{hard} constraint or a
\textit{contingent} feature of how mathematics has historically been taught. The response of
Section~8 is designed to be revisable either way: it protects construction where the dependency
binds, while allowing evidence to change the amount and form of unaided work required.

\subsection{Which Struggle Counts, and Why the Developmental Ordering Is Plausible}

The bottleneck rests on cognitive load theory, which at first appears to cut the other way.
Sweller's \citetext{1988} founding result was that conventional problem solving can
\textit{interfere} with learning. The means--ends search a novice runs consumes working memory
without building the schemas that learning consists in, so his remedy was to \textit{reduce}
problem-solving effort for novices, through goal-free problems and worked examples. Read
carelessly, that is an argument against the construction this paper insists on. Read carefully, it
is the argument's foundation, once Sweller's two kinds of load are separated.

\textit{Extraneous} load is capacity spent on processing that builds no schemas: the undirected
effort of an unguided means--ends search. \textit{Germane} load is the effortful processing that
does build them. Cognitive load theory objects to the first and depends on the second
\citetext{Sweller, 1988; Kirschner et al., 2006}. So ``construct before you delegate'' is not a
call to maximise difficulty. It is a call to preserve \textit{germane} load, the generative
processing through which schemas form.

It would be too quick to say the cognitive-load and productive-struggle traditions agree. They do
not. In its strong form, cognitive load theory holds that novices learn best from guided
instruction and worked examples and worst from unguided problem solving \citetext{Kirschner et
al., 2006}. The productive-failure programme holds almost the reverse: unguided problem solving
that \textit{fails} can build deeper, more transferable understanding than guided instruction,
provided structure follows the struggle \citetext{Kapur, 2008}. The disagreement is live and
unresolved. This paper need not settle it, because the two rivals agree on the one point that
matters here: the learner must engage in the generative work needed for learning. Sweller's worked examples help only when
the student self-explains them and then generates independently \citetext{Chi et al., 1989; Renkl,
2014}; Kapur's failure is productive only because the student did the struggling. Neither camp
endorses handing a novice a finished solution they never process.

Unstructured AI use can create precisely this problem. It is tempting to call AI a worked example with
the learner's generation removed, but that undersells the tool: AI does not merely present a finished
solution; it explains each step, links the underlying theory, and can be asked to refine its answer,
which makes it interactive rather than inert. That richness cuts against learning when the output is
simply received. The worked-example effect never came from the richness of the example but from the
student's own generation and self-explanation, and a more complete, more fluent solution removes that
work more thoroughly than a textbook example did, since the learner need not even assemble the
explanation or seek out the theory. The same interactivity could instead be turned the other way, to
withhold the answer and prompt the student to attempt, explain, and correct, so that the tool elicits
the learner's generation rather than supplying its own. Which it does is a matter of design, not a
property of the tool. The concern is therefore not AI assistance as such, but assistance that
routinely substitutes for the learner's generative processing. On the productive-failure side, AI is the instant, total
provision of structure, pre-empting the very struggle Kapur found productive. Kapur's result
pushes the proxy problem to its limit: the groups that produced the \textit{worse} immediate
solutions \textit{learned more}, beating the well-supported groups even on the problem type those
groups had practised. The visible artefact was not a weak proxy for competence but an inverted
one. Two qualifications keep the analogy honest. Kapur's setting is collaborative, ill-structured
group work, not the individual competence at issue here, so the transfer is by analogy, not
identity. And his failure is productive only because consolidation \textit{follows} it, as in
the sequence of Section~8, struggle is followed by feedback and verification, not left to run
forever. (That the advantage of studying over generating reverses as competence grows, the
expertise-reversal effect \citetext{Kalyuga et al., 2003}, is the same ordering from the other
side: study-heavy early, generation-heavy late.)

Cognitive load theory also supplies a mechanism that the taxonomies of Section~2 only describe, and
it puts that mechanism in a developmental relationship between kinds of activity rather than in a
ranking by difficulty. Bloom's hierarchy, as Section~2 noted, fits mathematics poorly, where an
``application'' may be trivial or research-level. What survives that critique is narrower: learners
need sufficiently automated and well-connected schemas for some forms of higher-order activity to be
manageable without excessive working-memory demands. This does not establish that unaided
hand-execution is the only route to those schemas, nor that more practice is always better. It does
establish a reason to be cautious about delegating the underlying work before students have acquired
the representations needed to understand and check it. And it sharpens the proxy problem. When AI performs the
Group~A operations, the student's schemas never automate, so the working-memory headroom for real
Group~C work never forms, yet the student can still produce artefacts that present as
higher-order: an analysis, an evaluation, a finished proof, because the machine supplies them. This
is the most misleading case of the proxy: work that looks like judgment while resting on none of the
construction it presupposes. Withdraw the AI and the apparent judgment collapses, because there was
never a secure procedural base beneath it.

\section{The Verification Instability}

\subsection{Oversight Requires Competence}

Let us consider the optimistic future approach, where AI generates, and the human verifies. Even if judgment stays with the person and execution goes to the machine, there is still a structural instability at its core. The issue is that verification is not one thing, and the division of labour is safe for only one kind of it.

Some checks are \textit{bounded}. They confirm a candidate answer cheaply and locally, without
reconstructing it: substitute a claimed eigenvector back into $Av=\lambda v$; substitute a candidate
solution back into an initial-value problem to see that the differential equation and the initial
data are satisfied; check that a probability lies in $[0,1]$; test a limiting case or an order of
magnitude. Where such a check exists, oversight really is cheaper than generation. A student who
could not have produced the result can still confirm it, and the optimistic division of labour
holds.

Other checks are \textit{unbounded}. Here there is no equation to substitute the answer into to see
whether it balances. Deciding that a model is the right model (whether a growth process is better
captured by an exponential or a logistic law), or that a boundary condition was appropriate,
or that a proof step is valid rather than merely plausible: none of these can be done without
essentially the competence needed to have done the work, and often more, because a plausible error
is harder to catch than to avoid. The tasks AI is most useful for (ill-structured modelling,
unfamiliar problems, research-level argument) are exactly the ones with no bounded check, and
there verification collapses back into full competence. So the oversight role is not uniformly
low-competence: it is low-competence where a cheap check happens to exist, and high-competence
precisely where AI is most tempting to trust and its errors least legible.

None of this is new outside mathematics. Bainbridge's \citetext{1983} ``ironies of automation'' made
the general point: automate the routine parts of a task and the operator is left responsible for what
cannot be automated, monitoring included, while disuse erodes the very skills monitoring needs, so the
better the automation, the less able the operator is to take over when it fails. Four decades of
human-factors work on automation complacency and skill decay confirm this across aviation, medicine,
and process control \citetext{Parasuraman \& Riley, 1997; Parasuraman \& Manzey, 2010}. Mathematical
verification under AI is the same irony moved into cognition: the competence to oversee AI is the
competence AI-assisted study is least likely to build.

There is a second layer, and it is worse. The ironies of automation assume skill erodes through
disuse, but oversight of a fluent system fails for a prior reason: people judge their own
comprehension badly, and the least competent judge it worst. Studies of comprehension calibration
find only a weak correlation between students' confidence that they have understood technical
material and their measured understanding \citetext{Maki, 1998, as reviewed in Person, Graesser, \&
Jordan, 2025}. The classroom illusion is the same finding in miniature: the sincere ``yes, I follow''
about material the student cannot use. A verifier who cannot gauge their own understanding cannot
reliably catch a plausible error in another's work, still less in the fluent, confident output of the
machine. So the instability is not only that the skill to verify erodes; it is that the capacity to
know \textit{when} verification is needed is what wholesale offloading leaves undeveloped.

The mathematics-specific evidence sharpens the point. Ask undergraduates to decide whether a
purported proof proves a theorem, and their accuracy can be poor; they may rely on surface features,
format, symbol use, or apparent authority rather than logical structure \citetext{Selden \& Selden,
2003}. This is the Section~2 incapacity from the other side. A student who can state the
$\varepsilon$--$\delta$ definition but cannot use it may also struggle to validate an AI-generated
argument using it. Proof validation is therefore a useful example of unbounded verification. The
important claim is not that such verification can only be learned by constructing proofs unaided,
but that meaningful validation requires mathematical resources that must themselves be developed.

\subsection{The Instability}

Together, these produce the \textit{verification instability}.

The reasoning becomes circular if verification is treated as automatically available. The
``AI generates, humans verify'' settlement is justified by the assumption that humans can verify
what the machine produces; yet some forms of verification require substantial prior mathematical
competence, and that competence may be weakened if the corresponding construction is routinely
offloaded (Section~5). The resulting tension is therefore conditional rather than universal. Where
a cheap, bounded check exists, delegation may be quite safe. Where verification is open-ended, the
ability to verify can approach the competence required to construct the work. This is why the paper
treats verification as a gating variable rather than assuming that human oversight will always be
available at low cost.

The early evidence fits the worry, with the usual caution. A field experiment with nearly a
thousand secondary mathematics students found that unrestricted access to a GPT-4 tutor raised
performance during practice and \textit{lowered} it on a later unaided exam, relative to students
who never had the tutor, while a tutor constrained to withhold answers and give hints did not show the same negative effect \citetext{Bastani et al., 2025}. One study, at secondary level, over a short
span: it invites replication before it carries much weight. But it is a direct experimental
instance of the dissociation the argument predicts, and its second finding, that constraint
removes the harm, is already evidence for the design response below. The cognitive-offloading
literature shows the same shape in memory: easy external access degrades both retention and
learners' awareness of what they do not know \citetext{Sparrow et al., 2011; Risko \& Gilbert,
2016}.

\subsection{Verification as the Gating Variable}

The instability has a constructive side, and the rest of the paper builds on it. How far a student
can \textit{independently verify} AI output is the natural measure of how much they can
\textit{safely delegate}. A student who can verify a routine integration may delegate routine
integration. A student who cannot has no business delegating it, because they cannot tell when the
machine is wrong. Verification competence is the gating variable for safe delegation, the
principle made central in Section~8.

\section{Design, Not Destiny (and Where the Stakes Are Sharpest)}

Nothing in the technology settles which future arrives. The same tool tends toward re-centering in
a struggle-first, construction-protecting design, and toward atrophy in an AI-permissive,
answer-on-demand default. The outcome is an institutional design choice. The central practical
claim of this paper is that mathematics education must make that choice deliberately, rather than
drift into the permissive default. That instructional design measurably moves learning, and not
only in principle, is well established: a meta-analysis across undergraduate STEM finds active,
construction-oriented teaching raises performance and lowers failure relative to lecturing
\citetext{Freeman et al., 2014}, and multi-institution studies of inquiry-based, struggle-first
college mathematics report gains in performance and persistence over lecture-based defaults
\citetext{Laursen et al., 2014; Kogan \& Laursen, 2014}.

The stakes are not uniform. They are sharpest where resources, infrastructure, and access are most
constrained, because the margin for a botched transition is smallest there. In much of the world, and acutely in many African universities, AI access is
stratified by connectivity, device ownership, and cost, so offloading is itself unevenly
distributed within a single cohort. Assessment combines AI-exposable coursework with unaided
high-stakes examinations, so competence gaps stay hidden until they are decisive. And importing AI
systems built elsewhere raises the question of whose mathematics is taught and certified. A full
treatment of these infrastructural and epistemic-political dimensions is beyond the scope of this
paper; here they enter for one reason: they make the design choice urgent rather than academic. The
gap between the two futures is widest, and hardest to absorb, exactly where resources are scarcest.

Language sharpens the stakes further. A mathematics student in a multilingual university navigates
not one language but several at once. There is the \textit{language of instruction}, the
official medium of lectures, textbooks, and examinations, typically a former colonial language such
as English, French, or Portuguese. There is the \textit{language of learning}, the language, or
blend of languages, in which the student actually reasons, questions, and makes meaning, often a
home language used alongside the official medium in the translanguaging characteristic of
multilingual classrooms \citetext{Setati, 2005; Garc\'ia \& Wei, 2014}. There is the
\textit{language of mathematics} itself, the symbolic, formal register that is no one's first
language and must be acquired as deliberately as any natural one \citetext{Pimm, 1987}. And now
there is a fourth: the \textit{language in which the student prompts AI}.

These layers are not incidental; they mediate every level of the approach. The gap between the
language of instruction and the language of learning imposes extra cognitive load exactly where
foundational construction (Level~1 of the sequence in Section~8) is most fragile
\citetext{Sweller, 1988; Adler, 2001}. A student building the
substitution schema or the eigenvector concept in a language they do not think in must spend
working memory on translation that a first-language peer spends on mathematics. The language of
mathematics compounds this, since the symbolic register must be built alongside the concepts it
expresses. Language, in short, raises the cost of building the base, and raises it most for the
students already least advantaged.

The prompting language adds a distinctive, underappreciated twist to the verification instability
of Section~6. The language in which a student queries AI partly determines the reliability of what
returns. Systems trained overwhelmingly on dominant-language corpora degrade markedly outside them,
so a student who prompts in their language of learning may receive \textit{more} error-prone
mathematics, and these are the students whose verification competence is already most strained
by the instruction-language gap. The
alternative, prompting in the language of instruction for more reliable output, buys
reliability at the price of a further translation, from the machine's register back into the
language in which the student must ultimately understand and verify the result; and every
translation is itself a site of error. Either way, the language layers do not sit beside the
verification problem. They compound it. Deliberate design must therefore treat language as a
first-class variable, deciding for each level which language construction, verification, and
delegation occur in.

\section{Two Principles: Construct Before You Delegate, and Let Verification Set the Ceiling}

\subsection{A Response Stated as Principles, Not Another Framework}

Section~4 noted that the field already has many named models of AI-era cognition, and part of the
complaint was that each adds structure without moving the argument forward. The response here is
therefore stated as principles rather than as a framework, because the two do different work. A
principle is a rule: it states what must hold. A framework is the organised set of activities through
which a principle is put into practice, and the same principle can be implemented by different
frameworks in different settings. Because the principle is the invariant and a framework is only one
contingent way of realising it, naming a framework would tie the argument to a particular
implementation, one that dates as the tools and the context change, and would add to a literature
already crowded with such models without stating anything the principle does not. Stating the
principles says what any workable implementation must respect and leaves institutions free to build
the framework that fits their own conditions. The principles are the contribution; the illustrative
sequence below is one way to see them at work, not a taxonomy to adopt.

\textbf{Construct before you delegate.} This principle follows from the bottleneck of Section~5, which
is a restriction on the acquisition of competence: where construction has not been sufficiently
developed, or where delegation takes place too early, the relevant skills either fail to develop or
develop only over more iterations. So before a class of mathematical work is routinely delegated to
AI, students should establish the conceptual and procedural resources needed to understand and verify
that class of work. In many cases this will require effortful construction and unaided practice; the
principle does not assume that every component must be mastered manually before any AI assistance is
permitted. Delegation is therefore conditioned on prior competence, not treated as a substitute for
building it.

\textbf{Let verification set the ceiling.} This principle follows from the verification instability of
Section~6. Where the underlying skills are not sufficiently developed, the learner cannot check the
accuracy of a result or judge whether the approach was even suitable, so what a student cannot verify
is exactly what they cannot safely delegate. Operationally the ceiling is a rule specific to the
student and the task: a student who can verify a result, not merely substitute an answer to recover a
balanced equation but relate the result to the problem and interpret it, has internalised the concepts
and may delegate up to that level; whatever sits above it must still be constructed by hand. The
ceiling sits at different heights for different work. A well-structured computation may carry a
bounded check, so a student who can run it may let AI do the arithmetic; a modelling problem requires
judging whether the model fits at all, which no plug-back settles, so a student who can only
substitute has not reached the ceiling that task demands. The same student therefore has a different
ceiling for different tasks, and delegation is safe only up to it.

Two further commitments follow as corollaries, not separate principles. Assessment must target the
construct, understanding, judgment, verification, formulation, and not the executable
procedure that served as its cheap proxy (Section~3); where procedure is assessed, it is assessed
as evidence of construction, under conditions that make construction visible. And the sequencing of
construction and delegation is a matter of deliberate design, not student discretion or
institutional drift (Section~7). Neither adds a new idea; each applies the two principles, to
assessment and to curriculum.

\subsection{An Illustrative Sequence}

The two principles imply an illustrative ordering of AI's role across a course, from protected
construction toward delegated work once the student can oversee what is handed over. Figure~2 sets
out one such sequence in five levels of competence, with the role AI might play at each, and places
it beside two instantiations, the well-structured eigenvalue problem and the ill-structured
harvesting problem developed below. It is an illustration of the construct-before-delegate
hypothesis, not a validated developmental model or a new taxonomy: nothing in the argument depends
on there being exactly five levels, or on learning running through them in a straight line.

\begin{center}
\vspace{4pt}
\begin{tikzpicture}[every node/.style={font=\scriptsize}]
\matrix (M) [
  matrix of nodes,
  nodes={draw, line width=0.4pt, rounded corners=2pt, inner sep=3pt, align=left,
         text width=3.9cm, anchor=north},
  column 1/.style={nodes={draw=none, text width=2.7cm}},
  row 1/.style={nodes={draw=none, font=\bfseries\scriptsize, align=center}},
  column sep=2mm, row sep=4mm,
]
{
 {} & The principle (abstract) & Well-structured: eigenvalue & Ill-structured: fishery \\
 \textbf{Level 1, Foundational construction.} \textit{AI withheld.} & Build the executable base by unaided, effortful construction. & Characteristic polynomial; eigenpairs $\lambda=1,3$ with $(1,\mp1)^{\mathsf T}$. & Solve the logistic model; equilibria $P_\pm$; maximum yield $H_{\max}=rK/4$. \\
 \textbf{Level 2, Conceptual consolidation.} \textit{AI deferred.} & Grasp \emph{why}; connect representations; form schemas. & Eigenvectors as preserved directions; symmetric $\Rightarrow$ real, orthogonal. & Stable vs.\ unstable equilibria; collapse threshold below $P_-$. \\
 \textbf{Level 3, Judgment \& strategy.} \textit{AI as partner.} & Choose the approach; judge plausibility. & Choose method; check $\operatorname{tr}A=4$, $\det A=3$. & Choose the model \& harvest form; judge parameter sensitivity. \\
 \textbf{Level 4, Verification \& critique.} \textit{AI as object of critique, the gate.} & Detect errors in AI output; judge arguments; this licenses delegation. & Catch a mis-paired eigenvector via $Av=\lambda v$ (plug-back). & Reject a single confident $H_{\max}$; interrogate the knife-edge assumption. \\
 \textbf{Level 5, Formulation \& creation.} \textit{AI as delegated tool.} & Formulate and create; oversee AI as a delegated tool. & Formulate $x_{k+1}=Ax_k$; delegate diagonalisation; interpret the dominant mode. & Formulate from messy data; delegate computation; give a conditioned range. \\
};
\foreach \c in {2,3,4}{%
  \foreach \r in {2,3,4,5}{%
    \pgfmathtruncatemacro{\rn}{\r+1}%
    \draw[-{Latex[length=1.8mm]}, line width=0.4pt] (M-\r-\c.south) -- (M-\rn-\c.north);%
  }%
}
\end{tikzpicture}

\vspace{4pt}
{\small\textbf{Figure 2: The Illustrative Sequence and Its Two Instantiations.} The same five
competence levels (rows), shown for the abstract principle and for two worked problems of opposite
character. AI's role expands down each column from \textit{withheld} to \textit{delegated}, with
Level~4 the gate that licenses delegation. Offered as illustration, not as a fixed model.\par}
\vspace{4pt}
\end{center}

\subsection{How the Sequence Resolves the Tension}

The sequence applies the two principles: AI's role expands from protected construction toward
delegated work in step with the student's developing competence and ability to verify. At Level~1,
AI is withheld from the protected construction activity. By Level~5, AI can be used as a delegated
tool for routine sub-tasks while the human formulates, judges, interprets, and verifies. Level~4 is
the pivot. Because verification sets the ceiling, demonstrated ability to critique AI output provides
the evidence needed before a class of work is routinely delegated. The point is not that every
student must pass through five universal stages, but that delegation should be conditioned on the
competence required to understand and check the delegated work.

\subsection{The Principles in Application: An Eigenvalue Problem}

The principles are easiest to grasp through worked examples, and we choose two of different
character, matched to the two schools of Section~5.3: one \textit{well-structured} problem, with a
canonical method and a definite answer, where AI's harm is to remove generation; and one
\textit{ill-structured} problem, where the assumptions must be made and reasonable answers differ,
where AI's harm is to manufacture false certainty. We trace each through the illustrative sequence,
showing how the treatment of the problem, and the role given to AI, change as the student ascends.
The first is a well-structured eigenvalue problem: diagonalise the symmetric matrix
\[
  A \;=\; \begin{bmatrix} 2 & 1 \\ 1 & 2 \end{bmatrix},
\]
that is, find its eigenvalues and eigenvectors. The object on the page is fixed; the cognitive
demand and AI's licence move with the level.

\subsubsection*{Level 1, Foundational construction (AI withheld)}
The student computes, unaided, the characteristic polynomial and solves it, then finds the
eigenvectors:
\[
  \det(A-\lambda I) = (2-\lambda)^2 - 1 = \lambda^2 - 4\lambda + 3 = (\lambda-1)(\lambda-3)
  \;\Rightarrow\; \lambda_1 = 1,\ \lambda_2 = 3,
\]
with eigenvectors $v_1=(1,-1)^{\mathsf T}$ for $\lambda_1=1$ and $v_2=(1,1)^{\mathsf T}$ for
$\lambda_2=3$. AI is withheld; the row-reduction and root-finding schemas are built here.

\subsubsection*{Level 2, Conceptual consolidation (AI deferred)}
The student explains \textit{why}: that an eigenvector is a direction the map merely scales, that
a real symmetric matrix is guaranteed real eigenvalues and orthogonal eigenvectors (indeed
$v_1\!\cdot\! v_2 = 0$ here), and that diagonalisation re-expresses the map in the basis where it
acts by simple scaling. The student writes this account first; AI then critiques it, rather than
producing it.

\subsubsection*{Level 3, Judgment and strategy (AI as partner)}
For an unfamiliar matrix the student must decide the approach, characteristic polynomial,
inspection of special structure (symmetric? triangular? a known form?), or a similarity argument, and sanity-check any result against the invariants: here the trace gives
$\operatorname{tr}A = 4 = \lambda_1+\lambda_2$ and the determinant gives
$\det A = 3 = \lambda_1\lambda_2$, both consistent. AI enters as a partner only after the student
has proposed a strategy and formed these expectations.

\subsubsection*{Level 4, Verification and critique (AI as object of critique)}
This is the gate. Suppose the AI returns the correct eigenvalues but pairs $\lambda=3$ with the
\textit{wrong} eigenvector, reporting
\[
  \lambda = 3,\quad v \overset{?}{=} (1,-1)^{\mathsf T} \qquad (\textbf{incorrect}).
\]
A student with secure Level~1--2 competence verifies against the defining relation $Av=\lambda v$:
\[
  A\begin{bmatrix}1\\-1\end{bmatrix} = \begin{bmatrix}1\\-1\end{bmatrix}
  = 1\cdot\begin{bmatrix}1\\-1\end{bmatrix} \;\neq\; 3\begin{bmatrix}1\\-1\end{bmatrix},
\]
so the reported pair fails: $(1,-1)^{\mathsf T}$ is the eigenvector for $\lambda=1$, not
$\lambda=3$. The correct pairing, $A(1,1)^{\mathsf T} = (3,3)^{\mathsf T} = 3(1,1)^{\mathsf T}$,
is confirmed the same way. A student who could not construct the eigenvectors cannot run this
check and will accept the plausible but transposed answer. Verification competence is what is
assessed, and what gates the next level.

\subsubsection*{Level 5, Formulation, modelling and creation (AI as delegated tool)}
The student is given a situation to model: two mutually reinforcing quantities evolve in discrete
steps, each period's values being $x_{k+1}=Ax_k$ for the coupling above, and the long-run
behaviour is required. The student must \textit{formulate} the model, identifying the state
vector, building the coupling matrix $A$ from the described interactions, and stating the
assumptions, and only then delegate the eigendecomposition to AI. Diagonalising $A=PDP^{-1}$
gives the closed form $x_k = P D^{k} P^{-1} x_0$, from which the student interprets the result:
the dominant eigenvalue $\lambda_2=3$ sets the asymptotic growth factor per step, and its
eigenvector $(1,1)^{\mathsf T}$ shows the two quantities approach equal proportion. The human owns
the formulation, the assumptions, the verification of the machine's decomposition, and the
interpretation; AI supplies only the mechanical decomposition the student has already, at
Level~1, proven able to supply.

\subsection{A Second Application: An Ill-Structured Harvesting Problem}

The eigenvalue problem is well-structured: one method, one answer, a verification that is a
plug-back check. The second example is its opposite, an \textit{ill-structured} modelling
problem in the sense of Section~5.3, where the model must be chosen, the assumptions supplied, and
reasonable answers differ. The setting is a fishery. A lake supports a fish stock that is being
harvested, and the question is how large a catch the stock can sustain. A natural first model adds
a harvest term to logistic growth,
\[
  \frac{dP}{dt} = rP\!\left(1-\frac{P}{K}\right) - H,
\]
with $P(t)$ the stock, $r$ and $K$ the intrinsic growth rate and carrying capacity, and $H$ the
constant catch per unit time. Unlike the eigenvalue problem, this one has no single right answer:
what counts as ``sustainable'' depends on assumptions the solver must make.

\subsubsection*{Level 1, Foundational construction (AI withheld)}
The student must do the underlying mathematics unaided. The unharvested logistic model is solved by
separation of variables and partial fractions, the integration that anchors the base, giving
$P(t)=K/(1+Ae^{-rt})$; and the harvested model's equilibria are found by setting $\tfrac{dP}{dt}=0$,
i.e.\ solving the quadratic $P^{2} - KP + KH/r = 0$, whose roots
\[
  P_\pm = \frac{K}{2}\left(1 \pm \sqrt{1 - \tfrac{4H}{rK}}\,\right)
\]
are real only when $H \le rK/4$. The maximum sustainable yield is therefore $H_{\max}=rK/4$, taken
at $P=K/2$, where the two equilibria coincide. AI is withheld; this analysis is the base.

\subsubsection*{Level 2, Conceptual consolidation (AI deferred)}
The student explains \textit{why}: for $H<H_{\max}$ the harvest cuts the growth curve at two points, an upper equilibrium $P_+$ that is stable and a lower one $P_-$ that is \textit{unstable}, a
collapse threshold below which the stock falls to zero. Harvesting at exactly $H_{\max}$ balances
the two on a knife-edge at $P=K/2$; harvesting above it drives the stock extinct. The student writes
this account first; AI then critiques it.

\subsubsection*{Level 3, Judgment and strategy (AI as partner)}
Here judgment dominates, and it is judgment about the \textit{model}, not the method. Is logistic
growth even appropriate, or does the stock need an age structure or a fluctuating environment? Is a
constant catch $H$ the right harvest form, or an effort-proportional one $hP$? Are the estimates of
$r$ and $K$, drawn from noisy, contested survey data, defensible, and how sensitive is the
answer to them? These are choices with no canonical resolution; the student forms a position and
brings it to AI as a partner, not a decision to be outsourced.

\subsubsection*{Level 4, Verification and critique (AI as object of critique)}
This is the gate, and its character has changed. Asked ``what catch is sustainable?'', AI will
readily return a single confident figure, typically $H_{\max}=rK/4$ for some assumed $r$ and $K$,
with the stock declared safe. There is no plug-back check that exposes the error, because the
arithmetic may be correct; the failure is in what the answer conceals. A student with secure
Level~1--2 competence interrogates the assumptions instead: the confident $H_{\max}$ sits exactly at
the knife-edge where the stable and unstable equilibria meet, so any shock, or any over-estimate of
$r$ or $K$, tips the stock past $P_-$ and collapses it. A defensible answer is a \textit{range}
below $H_{\max}$, with a precautionary margin, conditioned on the stated assumptions, not the
single number AI supplied. Verification here is not a substitution but a judgment about the
reasonableness of the model and its assumptions, and only a student who built the Level~1 analysis
can see that the single figure rests on a collapse threshold the model does not report.

\subsubsection*{Level 5, Formulation, modelling and creation (AI as delegated tool)}
Formulation \textit{is} the task. The student chooses the model and the harvest form, states the
assumptions (a closed population, a roughly constant environment, deterministic dynamics, a catch
that does not itself respond to the stock), estimates $r$ and $K$ from the available data together
with their uncertainty, and only then delegates the mechanical parts, solving the quadratic,
plotting the equilibria, running a sensitivity check, to AI. The human owns the model choice, the
assumptions, the interpretation, and, crucially, the honest representation of the answer as a
conditioned range rather than a point. What is assessed is formulation, assumption, oversight, and
the judgment to resist false precision, none of which the machine can supply.

\subsubsection*{What the two examples establish}
The pair spans the two schools of Section~5.3. In the well-structured eigenvalue problem, the key
distinction is between generation and a bounded verification check: the student can test
$Av=\lambda v$ without reconstructing every preceding step. In the ill-structured fishery problem,
by contrast, the mathematics may be correct while the model or its assumptions are inappropriate.
Verification therefore becomes a judgment about assumptions, sensitivity, and the meaning of the
result. The fishery example is consequently the sharper case for the paper's central claim:
\textit{correct mathematics is not necessarily correct modelling}. The same design principle applies
to both, but the competence required for verification rises with the openness of the task.

\section{Implications for Curriculum and Assessment}

\subsection{Curriculum}

A curriculum aligned to these principles specifies, for each competence, not only what students
should learn and how well, but \textit{when AI enters} and \textit{what prior competence makes its
use intelligible and checkable}. Early modules protect foundational construction and provide
opportunities for unaided attempts. Middle modules build conceptual consolidation and judgment, with
AI as deferred feedback and then dialogic partner. Later modules teach verification and critique,
including critique of AI output as a curricular objective, before routinely delegating classes of
work whose outputs students are expected to oversee. The curriculum therefore treats delegation as a
capability to be developed, not as a default condition of study.

\subsection{Assessment}

Assessment follows the construct-not-proxy corollary: it targets the construct, not the proxy. This
implies a bifurcation already visible across the sector, unaided, often oral or invigilated
assessment to certify individually held competence (especially the foundational construction of
Level~1 and the verification competence of Level~4), alongside AI-integrated, audit-trailed
assessment to develop and evidence the partnered and delegated competences of Levels~3 and~5. The
detailed design, worked templates and the equity constraints that shape them, and the
reconception of academic integrity it requires, from policing tool use to certifying authentic
cognitive work, lies beyond the scope of this paper. The point here is structural. Once
assessment targets the construct rather than the proxy, ``did the student use AI?'' loses its
centrality, because a well-designed assessment of this kind is valid whether they did or not.

\section{Implementing the Withheld Level When AI Is Everywhere}

The most immediate objection is practical and pointed. The first level demands that AI be
\textit{withheld} while foundational competence is built, yet AI now sits on every student's phone,
available at every moment outside the classroom. How can a competence be withheld when the tool is
ubiquitous? If this level cannot be realised, the approach collapses, because every level above it
depends on the base the withheld level builds. The objection deserves to be taken seriously.
Withholding is achievable, but only once ``withheld'' is understood correctly.

\subsection{What ``Withheld'' Cannot Mean}

Withholding cannot mean prohibition enforced by detection. That is the paradigm now widely shown to
be failing: detectors are unreliable, biased against second-language writers, and largely inapplicable
to symbolic mathematics; and prohibition merely drives use underground while punishing the visible
and the poorly resourced. An approach that depended on preventing students from touching AI in
their own time would be as unenforceable as the prohibition regimes it means to replace. If
``withheld'' meant ``technically prevented everywhere,'' the objection would be fatal. It does not.

\subsection{Withholding as Protected Construction Space, Not Surveillance}

Withholding is best understood spatially and temporally, not technologically: bounded, low-stakes,
in-class \textit{construction spaces} where the effortful building of foundational competence
actually happens, AI-free, because the activity is structured that way and the setting makes it
natural. Three design moves make this realistic.

\textit{First, relocate construction into contact time.} The single most consequential change is to
stop assigning foundational construction as unsupervised homework, where AI is unavoidable, and to move it into the lecture and tutorial, where it can be structured and observed. Board work,
live problem-solving on paper, think--pair--share, and worked-example practice done \textit{in the
room} convert construction from a private act the student may offload into a public one they
perform. The classroom becomes the site of doing, not merely of watching; the flipped emphasis is
not new \citetext{Freeman et al., 2014}, but AI gives it fresh urgency.

\textit{Second, decouple construction from grades.} Much offloading is a rational response to
high-stakes assessment: students offload because marks are at stake and the tool is available. If
foundational construction is practised in low-stakes, formative, ungraded or lightly-graded
settings, the incentive to offload largely dissolves, there is little to be gained by
outsourcing practice that is not being graded, while the stakes are relocated to the
AI-robust verification and oral checkpoints (Level~4) that cannot be faked. Withholding
works best where there is nothing to be won by defecting from it.

\textit{Third, make construction device-light by design.} Paper-and-board construction requires
neither connectivity nor a capable device, which is not only AI-free but an equity advantage: it
does not presuppose the infrastructure that is stratified across many African universities. In a resource-constrained setting, the low-technology construction space is both the
most AI-resistant and the most broadly accessible option, a case where the pedagogically and the
economically sound coincide.

\subsection{What This Means for the Lecture: A Migration of Function}

Relocating construction into contact time has an implication for the lecture, stated carefully. The
claim is not that lectures should disappear, but that generative AI changes the \textit{opportunity
cost} of contact time. The lecture historically fused transmission of content (exposition,
explanation, worked examples) with the modelling of expert thought, while the student's own
construction was displaced to private study. AI now makes explanation and worked examples easier to
obtain outside class, if unevenly. One design implication is to use more contact time for supervised
construction, diagnosis, and verification, where the educator can observe reasoning and respond to
misconceptions as they form, not an empirical finding that lectures have lost their value.

What the room then supplies is not merely active students but construction made \textit{visible to
the educator}. Where construction is homework, and may be AI-assisted, a competence gap stays invisible
until later assessment; a construction-mode room restores formative observation, letting the educator
watch a student get stuck and diagnose the misconception while it is still cheap to correct. The
room's comparative advantage is therefore as much \textit{diagnosis} as activity.

Two clarifications keep this from becoming either utopianism or a resourcing slogan. First,
construction mode is an activity design, not a particular staffing ratio; it can include peer
instruction and structured problem solving at lecture scale \citetext{Mazur, 1997}. Second, the
argument is compatible with continued exposition where exposition is the most effective way to
establish a prerequisite representation. The proposed migration is therefore selective: protect
contact time for forms of construction and diagnosis that are difficult to outsource, rather than
abolish transmission altogether.

\subsection{Recruiting Students into Withholding}

Even the best-designed construction space is porous; a determined student can still consult AI
before the tutorial. So the decisive lever is not enforcement but \textit{buy-in}: students
withhold AI from themselves once they understand why the construction matters. This is achievable,
and it is itself a teachable competence. Make the verification instability explicit, show
students the plausible, tidy, wrong AI solution of Section~8 and let them find that only the
competence they built themselves lets them catch it, and they have a first-hand reason to value
construction. So does \textit{metacognitive calibration}, the monitoring and self-assessment that
research on mathematical problem solving locates at the heart of competence \citetext{Schoenfeld, 1985}: brief, regular, unaided competence checks
that let students feel the gap between their assisted and unaided performance, countering the
illusion (Section~6) that fluency of access is possession of knowledge. A student who has felt that
gap has a motive to build the base that no rule can supply. Withholding, in the end, is something
students must be brought to choose; the educator's task is to make the choice intelligible, and the
classroom the easiest place to make it.

\subsection{The Honest Limit}

None of this guarantees that every student builds the base; no arrangement can, and the approach
should not pretend otherwise. What these moves achieve is more modest and more durable: they make
self-withholding the path of least resistance and greatest sense, they remove the incentives that
drive offloading, and, crucially, they relocate the \textit{certification} of foundational
competence to the AI-robust checkpoints of Level~4, where a student who skipped the construction
cannot conceal the fact. ``Withheld'' thus means not \textit{technically prevented} but
\textit{not required, not incentivised, structurally displaced into protected time, and ultimately
verified where it cannot be faked}. Understood this way, the withheld level is not a fantasy of
control in an AI-saturated world; it is a realistic design commitment that the rest of the
approach can stand on.

\section{Limitations and Conclusion}

\subsection{Limitations}

Three limitations warrant statement. The first is the most important. The approach rests on the
assumption that the construct-before-delegate dependency is substantially hard, and that assumption, as Section~5 stresses, is not proven. If future evidence shows that AI can scaffold robust
understanding along developmental routes that do not need the traditional volume of hand-execution,
the ordering the principles enforce would need revision. The approach is built to be revisable in
that event, but its present strictness is a bet on a still-open question. Second, the approach is a
theoretical proposal without empirical validation. The sequencing it prescribes must be tested
against learning outcomes, and that testing remains to be done. Third, the paper brackets the
infrastructural and epistemic-political dimensions that decide how the approach can actually be
enacted in resource-constrained and post-colonial settings. The design logic is general; its
implementation is not context-free.

\subsection{Conclusion}

Generative AI has not so much damaged university mathematics education as exposed a long-standing
weakness: the reliance on executable procedure as a convenient proxy for the understanding, judgment,
and verification the discipline actually values. Once AI can contribute to the observable performance,
the warrant for inferring individual competence from it must be reconsidered. The answer cannot be to
teach only the higher-order capacities hardest to delegate, because meaningful use of those capacities
may itself presuppose representations, procedures, and verification skills students still need to
develop. Yet the paper does not establish that traditional unaided execution is the only route to
them. That is an empirical question.

The two principles are therefore deliberately conditional: construct the conceptual and procedural
resources needed to understand a class of mathematical work before routinely delegating it, and let
independent verification set the ceiling on safe delegation. Assessment should target the construct
rather than the proxy, using unaided or controlled evidence where individually held competence must be
certified and AI-integrated work where partnered competence is the object. Whether the
construct-before-delegate dependency is a hard constraint or a contingent artefact of how mathematics
has been taught remains to be established; until that evidence is available, the paper argues for the
precautionary position that the competence needed to understand and verify a class of work should be
established before that work is routinely delegated.

\section*{Declarations}
\noindent\textbf{Funding.} The author received no specific funding for this work.\\[4pt]
\textbf{Competing interests.} The author declares no competing interests.\\[4pt]
\textbf{Data availability.} Not applicable: this is a conceptual paper and no datasets were generated or analysed.\\[4pt]
\textbf{Ethics approval.} Not applicable: the study involves no human participants or personal data.\\[4pt]
\textbf{Author contributions.} B.~Osano is the sole author and is responsible for the conception, argument, and writing of the manuscript.

\section*{Use of Generative AI}
The author used a generative AI assistant (Anthropic's Claude) during preparation of this
manuscript: for surfacing and organising relevant literature, and for sharpening the arguments. All AI-assisted text and reasoning were critically reviewed, revised, and, where necessary, reconstructed by the author, who verified every reference against primary sources and takes full responsibility for the content, arguments, and conclusions of the paper. The AI tool is not an author. No data were generated or analysed using AI.

\newpage
\section*{References}
\begin{reflist}
Adler, J. (2001). \textit{Teaching mathematics in multilingual classrooms}. Kluwer Academic Publishers.

AlphaProof \& AlphaGeometry teams. (2024). \textit{AI achieves silver-medal standard solving International Mathematical Olympiad problems} [Blog post]. Google DeepMind.\\ https://deepmind.google/blog/ai-solves-imo-problems-at-silver-medal-level/

Anderson, L. W., \& Krathwohl, D. R. (Eds.). (2001). \textit{A taxonomy for learning, teaching, and assessing: A revision of Bloom's taxonomy of educational objectives}. Longman.

Ayodele, K. P., Obayiuwana, E., Lawal, A. R., Bamimore, A., Offiong, F. B., \& Peter, E. A. (2026). \textit{Revising Bloom's taxonomy for dual-mode cognition in human--AI systems: The augmented cognition framework} [Preprint]. arXiv. https://doi.org/10.48550/arXiv.2602.00697

Bainbridge, L. (1983). Ironies of automation. \textit{Automatica}, \textit{19}(6), 775--779.

Bastani, H., Bastani, O., Sungu, A., Ge, H., Kabakc\i, \"O., \& Mariman, R. (2025). Generative AI without guardrails can harm learning: Evidence from high school mathematics. \textit{Proceedings of the National Academy of Sciences}, \textit{122}(26), e2422633122. \\https://doi.org/10.1073/pnas.2422633122

Bergqvist, E. (2007). Types of reasoning required in university exams in mathematics. \textit{The Journal of Mathematical Behavior}, \textit{26}(4), 348--370.

Biggs, J. (1996). Enhancing teaching through constructive alignment. \textit{Higher Education}, \textit{32}(3), 347--364.

Bjork, E. L., \& Bjork, R. A. (2011). Making things hard on yourself, but in a good way: Creating desirable difficulties to enhance learning. In M. A. Gernsbacher, R. W. Pew, L. M. Hough, \& J. R. Pomerantz (Eds.), \textit{Psychology and the real world} (pp. 56--64). Worth Publishers.

Boud, D., \& Falchikov, N. (2006). Aligning assessment with long-term learning. \textit{Assessment \& Evaluation in Higher Education}, \textit{31}(4), 399--413.

Bu\c{c}inca, Z., Malaya, M. B., \& Gajos, K. Z. (2021). To trust or to think: Cognitive forcing functions can reduce overreliance on AI in AI-assisted decision-making. \textit{Proceedings of the ACM on Human--Computer Interaction}, \textit{5}(CSCW1), Article 188. https://doi.org/10.1145/3449287

Chi, M. T. H., Bassok, M., Lewis, M. W., Reimann, P., \& Glaser, R. (1989). Self-explanations: How students study and use examples in learning to solve problems. \textit{Cognitive Science}, \textit{13}(2), 145--182.

Clark, L. (2025). \textit{Reframing Bloom's for the age of AI: A white paper for future-ready educators} [White paper]. Anthology. https://www.anthology.com/paper/reframing-blooms-for-the-age-of-ai-a-white-paper-for-future-ready-educators

Curry, H. B. (1951). \textit{Outlines of a formalist philosophy of mathematics}. North-Holland.

Darlington, E. (2013). The use of Bloom's taxonomy in advanced mathematics questions. In C. Smith (Ed.), \textit{Proceedings of the British Society for Research into Learning Mathematics}, \textit{33}(1).

Dubinsky, E. (1991). Reflective abstraction in advanced mathematical thinking. In D. Tall (Ed.), \textit{Advanced mathematical thinking} (pp. 95--126). Kluwer Academic Publishers.

Ellington, A. J. (2003). A meta-analysis of the effects of calculators on students' achievement and attitude levels in precollege mathematics classes. \textit{Journal for Research in Mathematics Education}, \textit{34}(5), 433--463.

Ernest, P. (1989). The knowledge, beliefs and attitudes of the mathematics educator: A model. \textit{Journal of Education for Teaching}, \textit{15}(1), 13--33.

Ernest, P. (1991). \textit{The philosophy of mathematics education}. Falmer Press.

Ernest, P. (1998). \textit{Social constructivism as a philosophy of mathematics}. State University of New York Press.

Freeman, S., Eddy, S. L., McDonough, M., Smith, M. K., Okoroafor, N., Jordt, H., \& Wenderoth, M. P. (2014). Active learning increases student performance in science, engineering, and mathematics. \textit{Proceedings of the National Academy of Sciences}, \textit{111}(23), 8410--8415. \url{https://doi.org/10.1073/pnas.1319030111}

Frieder, S., Pinchetti, L., Griffiths, R.-R., Salvatori, T., Lukasiewicz, T., Petersen, P. C., \& Berner, J. (2023). Mathematical capabilities of ChatGPT. \textit{Advances in Neural Information Processing Systems}, \textit{36}. https://doi.org/10.48550/arXiv.2301.13867

Garc\'ia, O., \& Wei, L. (2014). \textit{Translanguaging: Language, bilingualism and education}. Palgrave Macmillan.

Gonsalves, C. (2026). Generative AI's impact on critical thinking: Revisiting Bloom's taxonomy. \textit{Journal of Marketing Education}, \textit{48}(1), 4--19. https://doi.org/10.1177/02734753241305980

Hiebert, J., \& Grouws, D. A. (2007). The effects of classroom mathematics teaching on students' learning. In F. K. Lester (Ed.), \textit{Second handbook of research on mathematics teaching and learning} (pp. 371--404). Information Age Publishing.

Hmoud, M., \& Shaqour, A. (2024). AIEd Bloom's taxonomy: A proposed model for enhancing educational efficiency and effectiveness in the artificial intelligence era. \textit{The International Journal of Technologies in Learning}, \textit{31}(2), 111--128. https://doi.org/10.18848/2327-0144/CGP/v31i02/111-128

Jain, J., \& Samuel, M. (2025). \textit{Bloom meets Gen AI: Reconceptualising Bloom's taxonomy in the era of co-piloted learning} [Preprint]. \\Preprints.org. https://www.preprints.org/manuscript/202501.0271

Ji, Z., Lee, N., Frieske, R., Yu, T., Su, D., Xu, Y., Ishii, E., Bang, Y. J., Madotto, A., \& Fung, P. (2023). Survey of hallucination in natural language generation. \textit{ACM Computing Surveys}, \textit{55}(12), Article 248. https://doi.org/10.1145/3571730

Kalyuga, S., Ayres, P., Chandler, P., \& Sweller, J. (2003). The expertise reversal effect. \textit{Educational Psychologist}, \textit{38}(1), 23--31.

Kapur, M. (2008). Productive failure. \textit{Cognition and Instruction}, \textit{26}(3), 379--424.

Kirschner, P. A., Sweller, J., \& Clark, R. E. (2006). Why minimal guidance during instruction does not work: An analysis of the failure of constructivist, discovery, problem-based, experiential, and inquiry-based teaching. \textit{Educational Psychologist}, \textit{41}(2), 75--86.

Lithner, J. (2008). A research framework for creative and imitative reasoning. \textit{Educational Studies in Mathematics}, \textit{67}(3), 255--276.

Kogan, M., \& Laursen, S. L. (2014). Assessing long-term effects of inquiry-based learning: A case study from college mathematics. \textit{Innovative Higher Education}, \textit{39}(3), 183--199. \url{https://doi.org/10.1007/s10755-013-9269-9}

Laursen, S. L., Hassi, M.-L., Kogan, M., \& Weston, T. J. (2014). Benefits for women and men of inquiry-based learning in college mathematics: A multi-institution study. \textit{Journal for Research in Mathematics Education}, \textit{45}(4), 406--418. \url{https://doi.org/10.5951/jresematheduc.45.4.0406}

Long, C., Dunne, T., \& De Kock, H. (2014). Mathematics, curriculum and assessment: The role of taxonomies in the quest for coherence. \textit{Pythagoras}, \textit{35}(2), Article 240. \url{https://doi.org/10.4102/pythagoras.v35i2.240}

Maki, R. H. (1998). Test predictions over text material. In D. J. Hacker, J. Dunlosky, \& A. C. Graesser (Eds.), \textit{Metacognition in educational theory and practice} (pp. 117--144). Lawrence Erlbaum Associates.

Mazur, E. (1997). \textit{Peer instruction: A user's manual}. Prentice Hall.

Messick, S. (1989). Validity. In R. L. Linn (Ed.), \textit{Educational measurement} (3rd ed., pp. 13--103). American Council on Education / Macmillan.

OpenAI. (2023). \textit{GPT-4 technical report} [Technical report].\\ arXiv. https://doi.org/10.48550/arXiv.2303.08774

Parasuraman, R., \& Manzey, D. H. (2010). Complacency and bias in human use of automation: An attentional integration. \textit{Human Factors}, \textit{52}(3), 381--410.

Parasuraman, R., \& Riley, V. (1997). Humans and automation: Use, misuse, disuse, abuse. \textit{Human Factors}, \textit{39}(2), 230--253.

Person, N. K., Graesser, A. C., \& Jordan, C. E. (2025). Question-asking in peer tutoring. In D. R. Sanford, R. Hodges, \& M. Steiner (Eds.), \textit{The handbook of peer tutoring}. Bloomsbury Academic.

Pimm, D. (1987). \textit{Speaking mathematically: Communication in mathematics classrooms}. Routledge \& Kegan Paul.

Radmehr, F., \& Drake, M. (2019). Revised Bloom's taxonomy and major theories and frameworks that influence the teaching, learning, and assessment of mathematics: A comparison. \textit{International Journal of Mathematical Education in Science and Technology}, \textit{50}(6), 895--920. \url{https://doi.org/10.1080/0020739X.2018.1549336}

Renkl, A. (2014). Toward an instructionally oriented theory of example-based learning. \textit{Cognitive Science}, \textit{38}(1), 1--37.

Risko, E. F., \& Gilbert, S. J. (2016). Cognitive offloading. \textit{Trends in Cognitive Sciences}, \textit{20}(9), 676--688.

Schoenfeld, A. H. (1985). \textit{Mathematical problem solving}. Academic Press.

Selden, A., \& Selden, J. (2003). Validations of proofs considered as texts: Can undergraduates tell whether an argument proves a theorem? \textit{Journal for Research in Mathematics Education}, \textit{34}(1), 4--36.

Setati, M. (2005). Teaching mathematics in a primary multilingual classroom. \textit{Journal for Research in Mathematics Education}, \textit{36}(5), 447--466.

Sfard, A. (1991). On the dual nature of mathematical conceptions: Reflections on processes and objects as different sides of the same coin. \textit{Educational Studies in Mathematics}, \textit{22}(1), 1--36.

Shapiro, S. (2000). \textit{Thinking about mathematics: The philosophy of mathematics}. Oxford University Press.

Smith, G. H., Wood, L. N., Coupland, M., Stephenson, B., Crawford, K., \& Ball, G. (1996). Constructing mathematical examinations to assess a range of knowledge and skills. \textit{International Journal of Mathematical Education in Science and Technology}, \textit{27}(1), 65--77.

Sparrow, B., Liu, J., \& Wegner, D. M. (2011). Google effects on memory: Cognitive consequences of having information at our fingertips. \textit{Science}, \textit{333}(6043), 776--778.

Sweller, J. (1988). Cognitive load during problem solving: Effects on learning. \textit{Cognitive Science}, \textit{12}(2), 257--285.

Tall, D., \& Vinner, S. (1981). Concept image and concept definition in mathematics with particular reference to limits and continuity. \textit{Educational Studies in Mathematics}, \textit{12}(2), 151--169.

Thompson, A. G. (1992). educators' beliefs and conceptions: A synthesis of the research. In D. A. Grouws (Ed.), \textit{Handbook of research on mathematics teaching and learning} (pp. 127--146). Macmillan.

Trinh, T. H., Wu, Y., Le, Q. V., He, H., \& Luong, T. (2024). Solving olympiad geometry without human demonstrations. \textit{Nature}, \textit{625}, 476--482.

Vygotsky, L. S. (1978). \textit{Mind in society: The development of higher psychological processes}. Harvard University Press.

Zakariya, Y. F. (2021). Self-efficacy between previous and current mathematics performance of undergraduate students: An instrumental variable approach to exposing a causal relationship. \textit{Frontiers in Psychology}, \textit{11}, Article 556607. \url{https://doi.org/10.3389/fpsyg.2020.556607}

\end{reflist}
\end{document}